\documentclass{siamltex}

\usepackage{amssymb}
\usepackage{amsmath}
\usepackage[final]{graphicx}
\usepackage{psfrag}
\usepackage{epsfig}
\usepackage{latexsym}
\usepackage{exscale}

\newcommand{\bN}{\mathbb{N}}
\newcommand{\bR}{\mathbb{R}}
\newcommand{\bC}{\mathbb{C}}

\newcommand{\veps}{\varepsilon}
\newcommand{\vfi}{\varphi}

\newcommand{\cL}{\mathcal{L}}
\newcommand{\cLinv}{\mathcal{L}^{-1}}
\newcommand{\real}{{\rm Re}\,}
\newcommand{\acosh}{{\rm arccosh}\,}

\newtheorem{nota}{Remark}

\title{On the implementation of exponential methods for semilinear parabolic equations}

\author{Mar\'ia L\'opez-Fern\'andez \footnotemark[2]}

\date{July 27, 2008}

\begin{document}

\maketitle

\renewcommand{\thefootnote}{\fnsymbol{footnote}}

\footnotetext[2]{Departamento de Matem\'aticas, Universidad
Aut\'onoma de Madrid, Madrid, Spain. ~E-mail: {\tt
maria.lopez@uam.es}. Supported by DGI-MCYT under projects MTM
2005-00714 and MTM 2007-63257, cofinanced by FEDER funds, and the
SIMUMAT project S-0505/ESP/0158 of the Council of Education of the
Regional Government of Madrid (Spain).\\
July 28, 2008.}

\begin{abstract}
The time integration of semilinear parabolic problems by exponential
methods of different kinds is considered. A new algorithm for the
implementation of these methods is proposed. The algorithm evaluates
the operators required by the exponential methods by means of a
quadrature formula that converges like $O(e^{-cK/\ln K})$, with $K$
the number of quadrature nodes. The algorithm allows also the
evaluation of the associated scalar mappings and in this case the
quadrature converges like $O(e^{-cK})$. The technique is based on
the numerical inversion of sectorial Laplace transforms. Several
numerical illustrations are provided to test the algorithm.

\end{abstract}



\section{Introduction}\label{sec_intro}
The good numerical results obtained from the application of
exponential methods to the time integration of stiff semilinear
problems, have motivated much interest on these kind of methods
during the last years, see for instance \cite{BerSkaWri,CalPal06,
CoMa02, HoOs05imp, HoOs05exp, KaTre05}. The problems under
consideration can be written in the abstract format
\begin{equation}\label{ivp}
u'(t)= Au(t) + f(t,u(t)), \quad u(0)=u_0, \quad 0\le t\le T,
\end{equation}
where $A$ is a linear operator representing the highest order
differential terms and $f$ is a lower order nonlinear operator. The
solution to the initial value problem \eqref{ivp} is then given by
the variation of constants formula and most of the exponential
methods considered in the literature are constructed from this
representation of the solution.

Let us consider for instance the family of multistep exponential
methods developed in \cite{CalPal06} and briefly reviewed in
Section~\ref{subsec_msm}. Given a stepsize $h>0$, $n\ge 0$ and
approximations $u_{n+j}\approx u(t_{n+j})$, $t_{n+j} = (n+j)h$,
$0\le j\le k-1$, the $k$-step method approximates the solution $u$
of \eqref{ivp} at $t_{n+k}=(n+k)h$ by
\begin{equation}\label{multistep_intro}
u_{n+k} = \phi_0(k,hA) u_n + h \sum_{j=0}^{k-1} \phi_{j+1}(k,hA)
\Delta^j f_n,
\end{equation}
where $f_n=f(t_n,u_n)$, $\Delta$ denotes the standard forward
difference operator, and for $\lambda \in \bC$ and $k\ge 1$,
\begin{equation}\label{phis_multistep}
\phi_0(k,\lambda) = e^{k\lambda}, \quad \phi_j(k,\lambda) = \int_0^k
e^{(k-\sigma)\lambda} \left(\begin{array}{c} \sigma\\j-1
\end{array} \right)\,d\sigma, \quad 1\le j\le k.
\end{equation}
As we can see in \eqref{multistep_intro}, these methods require the
evaluation of $\phi_j(k,hA)$, $0\le j\le k$, for $\phi_j(k,\lambda)$
defined in \eqref{phis_multistep}. This is in fact the main
difficulty in the implementation of the methods in
\eqref{multistep_intro} and, in general, of exponential methods,
since they typically require the evaluation of vector-valued
mappings $\phi(hA)$, with $h$ the time step in the discretization
and either
\begin{equation}\label{Phiexponencial}
\phi(\lambda) = e^{m\lambda}, \qquad \lambda \in \bC,
\end{equation}
or
\begin{equation}\label{Phigeneral}
\phi(\lambda) = \int_0^m e^{(m-\sigma)\lambda}p(\sigma)\, d\sigma,
\end{equation}
with $m$ an integer and $p(\sigma)$ a polynomial. The values of $m$
and $p$ in \eqref{Phigeneral} depend on the method. For instance, in
the case of the methods in \eqref{multistep_intro}, it is clear from
\eqref{phis_multistep} that $m=k$, the number of steps of the
method, and
$$
p(\sigma)=\left(\begin{array}{c} \sigma\\j-1
\end{array} \right) = \frac{\sigma(\sigma-1)\dots (\sigma-j+2)}{(j-1)!}, \qquad 1\le j\le k.
$$
For the explicit exponential Runge--Kutta methods constructed in
\cite{HoOs05exp} it is $m=1$ and
$$
p(\sigma) = \frac{\sigma^{k-1}}{(k-1)!}, \qquad k\ge 1,
$$
as we can see in Section~\ref{subsec_rkm}, \eqref{phis_rk}.

In the present paper we propose a way to evaluate the operators
$\phi(hA)$ in \eqref{Phiexponencial} and \eqref{Phigeneral} when $A$
in \eqref{ivp} is the infinitesimal generator of an analytic
semigroup in a Banach space $X$. Thus, we will assume that $A:
D(A)\subset X \to X$ is sectorial, i.e., $A$ is a densely defined
and closed linear operator on $X$ and there exist constants $M>0$,
$\gamma \in \bR$, and an angle $0 < \delta < \frac \pi 2$, such that
the resolvent fulfils
\begin{equation}\label{sectorial}
\| (zI - A)^{-1} \| \le \frac{M }{ |z-\gamma|}, \quad\ \hbox{ for
}\quad |\arg(z -\gamma)| < \pi-\delta.
\end{equation}

As a side product we also obtain an accurate way to evaluate the
mappings $\phi(\lambda)$ in \eqref{Phigeneral} at the scalar level,
which is itself a well-known problem in numerical analysis. This is
exemplified in \cite{KaTre05} with the mapping
\begin{equation}\label{vfi_1} \vfi(\lambda) =
\frac{e^{\lambda}-1}{\lambda},
\end{equation}
which is required for instance by the exponential Runge--Kutta
methods of \cite{HoOs05exp}. The evaluation of $\vfi$ for small
$\lambda$ by using formula \eqref{vfi_1} suffers from cancellation
error. On the other hand, the use of a truncated Taylor expansion
only works well for small enough $\lambda$. This implies that for
some intermediate values of $\lambda$ it is not very much clear how
to choose the proper formula and moreover both of them could lose
accuracy. For $\lambda$ inside a sector $|\arg(\lambda-\gamma)|\ge
\pi-\delta$, these difficulties can be overcome by writing $\vfi$ in
the format \eqref{Phigeneral}, with $m=1$ and $p(\sigma)=1$. By
doing so, we will be able to evaluate $\vfi(\lambda)$, independently
of the size of $\lambda$, by using essentially the same technique
developed in principle to evaluate the vector-valued mapping
$\vfi(hA)$.

In the recent literature, several alternatives have been proposed to
evaluate the required operators $\phi(hA)$ or alternatively their
action on given vector, based on a Krylov approach
\cite{HoLu97,HoLu98}, on Pad\'e approximants \cite{BerSkaWri} or on
a suitable contour integral representation of the mappings by means
of the Cauchy integral formula \cite{KaTre05}. In particular, in
\cite{KaTre05}, the goal is both the evaluation of the mappings
$\phi(h\lambda)$ at the scalar level, assuming that a
diagonalization of $A$ is available, but also at the operator level,
since it also allows the evaluation of the operators $\phi(hA)$
working with the full matrix $A$. However, despite the good
computational results reported in \cite{KaTre05}, the way of
selecting the parameters involved in the quadrature formulas is not
very much clear and they depend very strongly on the equation
considered and the spatial discretization parameters. The algorithm
we propose is derived by using Laplace transformation formulas and
is finally based on a suitable contour integral representation of
the mappings $\phi(hA)$, too. However, the quadrature formulas we
use borrow their parameters from the method in \cite{LoPSch} for the
numerical inversion of sectorial Laplace transform, where a rigorous
analysis of the error is performed together with an optimization
process to choose the different parameters involved in the
approximation.

In order to apply the quadrature formulas developed in
\cite{LoPSch}, we derive a representation of the operators in
\eqref{Phiexponencial} and \eqref{Phigeneral} as the inverses of
suitable Laplace transforms $\Phi(z,hA)$ at certain values of the
original variable $\sigma$. These Laplace transforms have all the
form
\begin{equation}\label{LTPhi_general}
\Phi(z,hA) = R(z)(zI-hA)^{-1},
\end{equation}
with $R(z)$ a scalar rational mapping of $z\in \bC$. Due to
\eqref{sectorial}, the mappings $\Phi(z,hA)$ turn out to be {\em
sectorial} in the variable $z$, i.e., there exist constants $\gamma
\in \bR$ and $M > 0$, possibly different from the constants in
\eqref{sectorial}, such that
\begin{equation}\label{LTsectorial}
  \begin{array}{c}
    \hbox{$\Phi(z,hA)$ is analytic for $z$ in the sector $|\arg(z-\gamma)| < \pi-\delta$
      and there}
    \\[2mm]\displaystyle
    \| \Phi(z,hA)\| \le \frac {M}{|z -
\gamma|^{\nu}}, \qquad\hbox{for some } \nu\ge 1.
  \end{array}
\end{equation}
In this way, we reduce the problem of computing $\phi(hA)$ to the
inversion of a sectorial Laplace transform $\Phi(z,hA)$ of the form
\eqref{LTPhi_general}. We then use the method for the numerical
inversion of sectorial Laplace transforms developed in \cite{LoPSch}
and reviewed in Section~\ref{sec_LT}. The inversion method consists
on a quadrature formula to discretize the inversion formula
\eqref{invLT}, so that we finally approximate
\begin{equation}\label{quadgeneral}
\phi(hA) \approx \sum_{\ell=-K}^{K} w_{\ell}e^{mz_{\ell}}
\Phi(z_{\ell},hA) = \sum_{\ell=-K}^{K} w_{\ell}e^{mz_{\ell}}
R(z_{\ell})(z_{\ell}I-hA)^{-1},
\end{equation}
with the quadrature weights $w_{\ell}$ and nodes $z_{\ell}$ given in
\eqref{pesos_nodos}. The convergence results in \cite{LoPSch} assure
an error estimate in the approximation \eqref{quadgeneral} at least
like $O(e^{-cK/\ln K})$. Further, by following this approach, our
selection of parameters in the implementation of exponential methods
depends only on $\delta$ in \eqref{sectorial}, being independent of
$h$ and $M$. In the particular case that we want to evaluate $\phi$
at a scalar $\lambda$ in the sector $|\arg(\lambda-\gamma)|\ge
\pi-\delta$, the approximation becomes simply
\begin{equation}\label{quadscalar}
\phi(\lambda) \approx \sum_{\ell=-K}^{K} w_{\ell}e^{mz_{\ell}}
\Phi(z_{\ell},\lambda) = \sum_{\ell=-K}^{K} w_{\ell}e^{mz_{\ell}}
\frac{R(z_{\ell})}{z_{\ell}-\lambda}
\end{equation}
and the convergence rate improves to an $O(e^{-cK})$.

The approximation in \eqref{quadgeneral} allows the computation of
all the operators $\phi(hA)$ required by an exponential method
before the time-stepping begins, so that only the products
matrix$\times$vector need to be implemented at every time step.
However, formula \eqref{quadgeneral} requires the computation and
storage of all the inverses $(z_{\ell}I-hA)^{-1}$, $\ell =
-K,\dots,K$. Even if $A$ is a sparse matrix and these inversions can
be performed efficiently, the storage of the resulting full matrices
$\phi(hA)$ and the subsequent products matrix$\times$vector can
become prohibitive for large problems. In this situation,
\eqref{quadgeneral} should be combined with a data sparse procedure
to approximate the resolvent operators, as it is proposed in
\cite{GavHaKh04,GavHaKh05}.

Another way to avoid the computation and storage of all the
resolvents in \eqref{quadgeneral} could be the application of the
formula to evaluate the action $\phi(hA)v$ on a given vector $v$,
instead of the full operator $\phi(hA)$. In this case, it is not
necessary the computation of the full inverses
$(z_{\ell}I-hA)^{-1}$, but the resolution of the linear systems
$$
(z_{\ell}I-hA)x = v.
$$
However, solving all these linear systems for all the quadrature
nodes $z_{\ell}$ in \eqref{quadgeneral} at every time step can
become quite expensive and, in this sense, the Krylov approach
developed in \cite{HoLu97} appears as a more competitive
alternative. In this situation, only formula \eqref{quadscalar} can
be useful, in order to evaluate the mappings $\phi(h\lambda)$ at the
eigenvalues of the Hessenberg matrices involved in the Krylov
approximation.

By using \eqref{quadgeneral} we are in fact computing an
approximation to the numerical solution of \eqref{ivp} provided by
an exponential method. Thus, the global error after applying
\eqref{quadgeneral} to the time integration of \eqref{ivp} can be
split into the error in the time integration by the ``pure''
exponential method and the deviation from the numerical solution
introduced by the approximation \eqref{quadgeneral} of the operators
$\phi(hA)$. The error in the time integration for the exponential
integrators considered in the present paper is analyzed in
\cite{CalPal06} and \cite{HoOs05exp} (see
Section~\ref{sec_expmethods}), and the quadrature error is analyzed
in \cite{LoPSch} (see Section~\ref{sec_LT}, Theorem~\ref{thm:err}).
In order to visualize the effect of this approximation, we show in
Section~\ref{sec_experiments} the performance of our implementation
for several problems with known exact solution and moderate size
after the spatial discretization. In the error plots provided we can
observe that the error coincides with the expected error for the
exact exponential integrators up to high accuracy for quite moderate
values of $K$ in \eqref{quadgeneral}, i.e., the error induced by the
quadrature \eqref{quadgeneral} is negligible compared to the error
in the time integration.

Finally we notice that the matrix exponential $e^{tA}$ and also
certain rational a\-ppro\-xi\-ma\-tions to it originating from
Runge--Kutta schemes have already been successfully approximated by
using this approach \cite{LoLuPaSch,LoP04,LoPSch}.

The paper is organized as follows. In Section~\ref{sec_expmethods}
we consider the class of explicit multistep exponential methods
proposed in \cite{CalPal06} and the explicit exponential
Runge--Kutta methods in \cite{HoOs05exp}. Section~\ref{sec_LT} is a
review of the method for the numerical inversion of sectorial
Laplace transforms presented in \cite{LoPSch}. In
Section~\ref{sec_eval} we deduce a representation for the operators
required in the implementation of these integrators in terms of
suitable Laplace transforms and apply the inversion method to the
implementation of exponential methods. In
Section~\ref{sec_evalscalar} we consider with some detail the
evaluation of the associated mappings at the scalar level and
present some numerical results. We finally test our algorithm with
several examples in Section~\ref{sec_experiments}, where we
implemented \eqref{quadgeneral} by using the full matrices.

\section{Exponential methods}\label{sec_expmethods}
In this section we review some of the exponential methods in the
recent literature. In this way, we consider the class of explicit
multistep exponential integrators developed in \cite{CalPal06} and
the exponential Runge--Kutta methods of \cite{HoOs05exp}, which are
explicit as well.

\subsection{Multistep exponential methods}\label{subsec_msm}
In \cite{CalPal06}, a class of explicit exponential methods is
constructed for the time integration of problems of the form
\eqref{ivp} with $A:D(A)\subset X \to X$ the infinitesimal generator
of a $C_0$-semigroup $e^{tA}$, $t\ge 0$, of linear and bounded
operators in a Banach space $X$.

As we already stated in the Introduction, we will restrict ourselves
to the case of $A$ in \eqref{ivp} sectorial. Then, for $0\le \alpha
<1$, the fractional powers $(\omega-A)^{\alpha}$ are well defined
for $\omega
> \gamma$ in \eqref{sectorial}, and $X_{\alpha} =
D((\omega-A)^{\alpha})$ endowed with the graph norm
$\|\cdot\|_{\alpha}$ is a Banach space \cite{Henry}. In this
situation, the nonlinearity $f$ in \eqref{ivp} is assumed to be
defined on $[0,T]\times X_{\alpha} \to X$, for some $0\le \alpha <
1$, and to be locally Lipschitz in a strip along the exact solution.
Thus, there exists $L(R,T)>0$ such that
\begin{equation}\label{Lipschitz}
\|f(t,\eta)-f(t,\xi)\| \le L \|\eta-\xi\|_{\alpha}, \quad \eta, \xi
\in X_{\alpha}, \quad 0\le t\le T,
\end{equation}
for $\max\left( \|\eta-u(t)\|_{\alpha}, \|\xi-u(t)\|_{\alpha}
\right) \le R$.

For $k\ge 1$, the $k$-step method is constructed from the variation
of constants formula in the interval $[t_n, t_{n+k}]$ as follows.
Taking a stepsize $h=T/N$, $N\ge k$, and the corresponding time
levels $t_n=nh$, $0\le n\le N$, the solution $u$ of \eqref{ivp} at
$t_{n+k}$ is given by
\begin{equation}\label{vcf}
u(t_{n+k}) = e^{khA}u(t_n) + h \int_0^k e^{(k-\sigma)hA}
f(t_n+\sigma h, u(t_n + \sigma h))\, d\sigma.
\end{equation}
Given approximations $u_{n+j}\approx u(t_{n+j})$, $0\le j \le k-1$,
the approximation $u_{n+k}\approx u(t_{n+k})$ is obtained after
replacing $f$ in \eqref{vcf} by the Lagrange interpolation
polynomial of degree $k-1$, $P_{n,k-1}$ through the points
$\{(t_{n+j},f(t_{n+j},u_{n+j}))\}_{j=0}^{k-1}$ and integrating
exactly. Writing
\begin{equation}
P_{n,k-1}(t_n + \sigma h) = \sum_{j=0}^{k-1} \left(\begin{array}{c}
\sigma\\j \end{array} \right) \Delta^{j}f_n,
\end{equation}
with $f_m=f(t_m,u_m)$, $0\le m\le N-1$, and $\Delta$ the standard
forward difference operator, the approximation $u_{n+k}$ is computed
from
\begin{equation}\label{multistep}
u_{n+k} = \phi_0(k,hA) u_n + h \sum_{j=0}^{k-1} \phi_{j+1}(k,hA)
\Delta^j f_n,
\end{equation}
where, for $\lambda \in \bC$, $k\ge 1$, and $0\le j\le k$, the
mappings $\phi_j(k,\lambda)$ are given in \eqref{phis_multistep}.
The methods defined in \eqref{multistep} are explicit and, as we
already mentioned in the Introduction, they require the evaluation
of $\phi_j(k,hA)$.

In \cite[Theorem 1]{CalPal06} it is shown that if the nonlinearity
$f(t,u(t))$ belongs to the space $C^k([0,T],X)$ and the starting
values $u_0,\dots,u_{k-1}\in X_{\alpha}$ satisfy
$$
\|u(t_j) - u_j\|_{\alpha} \le C_0 h^k, \qquad 0\le j\le k-1,
$$
the method defined in \eqref{multistep} exhibits full order $k$,
i.e., there exists $C > 0$ such that
\begin{equation}\label{errbound}
\| u(t_n) - u_n \|_{\alpha} \le C h^k \|f^{(k)}\|_{\infty}, \qquad
0\le n\le N.
\end{equation}
The constant $C$ depends on $k,\alpha,T$, $L(R,T)$ in
\eqref{Lipschitz}, and $\gamma, M$ in \eqref{sectorial}, but it is
independent of $h$ and $f$.

\subsection{Exponential Runge--Kutta methods}\label{subsec_rkm}
Explicit exponential Run\-ge--Ku\-tta methods are presented in
\cite{HoOs05exp} for the time integration of semilinear parabolic
problems. For $h=T/N$, $N \ge 1$, and $1\le i\le s$, the
approximations $u_n$ to $u(t_n)$, with $t_n=nh$, are given by
\begin{equation}\label{rk}
\begin{array}{rcl}
U_{ni} &=& \displaystyle e^{c_i hA}u_n + h
\sum_{j=1}^{i-1}a_{ij}(hA)f(t_n+c_j
h,U_{nj}), \\
u_{n+1} &=& \displaystyle e^{hA}u_n + h \sum_{i=1}^s
b_i(hA)f(t_n+c_i h,U_{ni}),
\end{array}
\end{equation}
with $c_1=0$ ($U_{n1}=u_n$). The coefficients $b_i(\lambda)$ and
$a_{ij}(\lambda)$ are linear combinations of $\vfi_k(\lambda)$ and
$\vfi_k(c_l\lambda)$ with
\begin{equation}\label{phis_rk}
\vfi_k(\lambda) = \int_0^1 e^{(1-\sigma)\lambda}
\frac{\sigma^{k-1}}{(k-1)!}\,d\sigma, \qquad \lambda\in \bC, \quad
k\ge 1, \quad t
> 0.
\end{equation}
Setting $\vfi_0(\lambda) = e^{\lambda}$, we see that the
implementation of \eqref{rk} requires the e\-va\-lua\-tion of
$\vfi_k(hA)$ and $\vfi_k(c_lhA)$, for $1\le l \le s$ and several
values of $k\ge 0$.

As in \cite{CalPal06}, the nonlinearity $f$ in \eqref{ivp} is
assumed to satisfy a local Lipschitz condition like
\eqref{Lipschitz} and the solution $u:[0,T]\to X_{\alpha}$ and
$f:[0,T]\times X_{\alpha} \to X$ are assumed to be sufficiently
smooth so that the composition
$$
f^*:[0,T]\times X_{\alpha}\to X \, :\, t\to f^*(t)=f(t,u(t))
$$
is a sufficiently often differentiable mapping, too. Under these
assumptions, stiff order conditions are derived and exponential
Runge--Kutta methods of the form \eqref{rk} are constructed up to
order four in \cite{HoOs05exp}.

\section{The numerical inversion of sectorial Laplace
transforms}\label{sec_LT} In this section we review the numerical
inversion method for sectorial Laplace transforms presented in
\cite{LoPSch}, which is a further development of \cite{LoP04}.

\begin{figure}
\centering
\includegraphics[width=.6\textwidth]{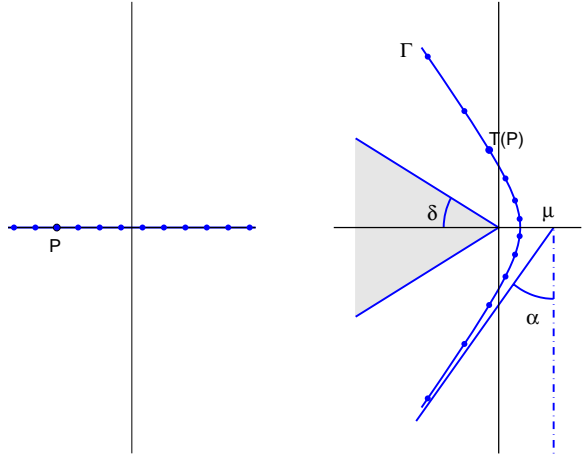}
\caption{Action of $T$ in \eqref{hyppara} on the real
axis}\label{fig:hyp}
\end{figure}

For a locally integrable mapping $f:(0,\infty) \to X$, bounded by
\begin{equation}\label{expgrowtn}
\|f(t)\| \le C t^{\nu-1} e^{\gamma t}, \qquad \mbox{for some } \
\gamma\in \bR, \ \nu >0,
\end{equation}
we denote its Laplace transform
\begin{equation}\label{LT}
F(z) = \cL[f](z) = \int_0^{\infty} e^{-tz}f(t)\, dt, \qquad \real z
> \gamma.
\end{equation}
When $F$ satisfies \eqref{LTsectorial}, the method in \cite{LoPSch}
allows to approximate the values of $f$ from few evaluations of $F$.
This is achieved by means of a suitable quadrature rule to
discretize the inversion formula
\begin{equation}\label{invLT}
f(t) = \frac 1{2\pi i} \int_{\Gamma} e^{zt} F(z)\,dz,
\end{equation}
where $\Gamma$ is a contour in the complex plane, running from
$-i\infty$ to $i\infty$ and laying in the analyticity region of $F$.
Due to \eqref{LTsectorial}, $\Gamma$ can be taken so that it begins
and ends in the half plane $\real z <0$. Following \cite{LoPSch}, in
\eqref{invLT} we choose $\Gamma$ as the left branch of a hyperbola
parameterized by
\begin{eqnarray}\label{hyppara}
\bR \to \Gamma &:& \ x \mapsto T(x) =
  \mu (1 - \sin(\alpha+ix)) + \gamma,
\end{eqnarray}
where $\mu > 0$ is a scale parameter, $\gamma$ is the shift in
\eqref{LTsectorial}, and $0 < \alpha < \frac{\pi}2 -\delta$. Thus,
$\Gamma$ is the left branch of the hyperbola with center at
$(\mu,0)$, foci at $(0,0)$, $(2\mu,0)$, and with asymptotes forming
angles $\pm(\pi/2 + \alpha)$ with the real axis, so that $\Gamma$
remains in the sector of analyticity of $F$,
$|\arg(z-\gamma)|<\pi-\delta$. In Figure~\ref{fig:hyp} we show the
action of the conformal mapping $T$ on the real axis.

After parameterizing \eqref{invLT}, the function $f$ is approximated
by applying the truncated trapezoidal rule to the resulting integral
along the real axis, i.e.,
\begin{eqnarray}\label{numinvLT}
f(t) = \frac{1}{2\pi i} \int_\Gamma   e^{tz}\,F(z)\, dz \approx
\sum_{\ell=-K}^K w_{\ell}\, e^{tz_{\ell}}\,F(z_{\ell}),
\end{eqnarray}
with quadrature weights $w_{\ell}$ and nodes $z_{\ell}$ given by
\begin{equation}\label{pesos_nodos}
w_{\ell} = -{\tau \over 2 \pi i}\: T'(\ell\tau)~, \quad z_{\ell} =
T(\ell \tau), \quad -K \le \ell \le K,
\end{equation}
and $\tau > 0$ a suitable step length parameter. We notice that the
minus sign in the formula for the weights comes from setting the
proper orientation in the parametrization of $\Gamma$. In case of
symmetry, the sum in \eqref{numinvLT} can be halved to
\begin{equation}\label{numinvLT_real}
f(t) \approx \, \real \left(\sum_{\ell = 0}^{K} w^*_{\ell}
e^{tz_{\ell}} F(z_{\ell}) \right),
\end{equation}
with $w^*_0 = w_0$ and $w^*_{\ell} = 2w^*_{\ell}$, $\ell \ge 1$. The
good behavior of the quadrature formula \eqref{numinvLT} is due to
the good properties of the trapezoidal rule when the integrand can
be analytically extended to a horizontal strip around the real axis
\cite{Ste,Sten}.

During the last years, different choices of contours $\Gamma$ and
parameterizations have been studied for the numerical inversion of
sectorial Laplace transforms. Apart from the approach in
\cite{LoP04,LoPSch}, which is the one we follow, the choice of a
hyperbola has been studied in
\cite{GavHaKh04,GavHaKh05,Mcl,She,WeiT}. The choice of $\Gamma$ as a
parabola has been considered recently in
\cite{GavHaKh04,GavHaKh05,WeiT} and we refer also to Talbot's method
\cite{Talbot79,Wei} for another kind of integration contour
$\Gamma$, with horizontal asymptotes as $|z|\to \infty$.

The error analysis for \eqref{numinvLT} shows that the role played
by the round-off errors is very important. Several ways to avoid
error propagation are proposed in \cite{LoPSch}, depending on the
available information about the errors in the evaluations of
$F(z_{\ell})$ and the elementary mappings involved in
\eqref{numinvLT}. Our algorithm for exponential methods uses
\eqref{quadgeneral} to approximate the required operators $\phi(hA)$
in \eqref{Phiexponencial} and \eqref{Phigeneral}, i.e.,
\eqref{numinvLT} with Laplace transforms of the form
\eqref{LTPhi_general}. Thus, the required Laplace transforms will be
in practice evaluated by performing the inversion of few matrices of
the form $zI-hA$, with $A$ being normally a discrete version of the
operator in \eqref{ivp} and $h$ the time step length. Taking into
account that this is likely to be accomplished by means of some
auxiliary routine, we propose to use in this context the least
accurate version of the method in \cite{LoPSch}, where the error
propagation is avoided without using any information about the
errors in the evaluations of \eqref{numinvLT}.

The following result provides an error bound for the proposed
version of the inversion method like $O(e^{-cK/\ln K})$, with $K$
the number of quadrature nodes.
\begin{theorem}\label{thm:err}\cite{LoPSch} Assume that the Laplace transform $F(z)$
satisfies the sectorial condition \eqref{LTsectorial} 
and let $\alpha$ and $d$ be such that
\begin{equation}\label{angulos}
0 < \alpha-d < \alpha+d < \frac{\pi}2-\delta.
\end{equation}
For $t_0 > 0$, $\Lambda \ge 1$ and $K\ge 1$, we select the
parameters
\begin{equation}\label{parameters}
\tau = \frac {a(K)}{K}, \qquad \mu = \frac{2\pi d}{\Lambda t_0
a(K)},
\end{equation}
with $a(K)= \acosh \left(\Lambda K/\sin \alpha \right)$.

Then, the error $E_K(t)$ in the approximation \eqref{numinvLT} to
$f(t)$ with quadrature weights and nodes in \eqref{pesos_nodos} is
bounded by
\begin{equation}\label{cota}
\|E_K(t)\| \le M\, \Pi\, Q\, e^{2\pi d/a(K)}\, t^{\nu-1} \left(
\veps + \frac{e^{-2\pi d K/a(K)}}{1-e^{-2\pi d K/a(K)}} \right),
\end{equation}
uniformly for $t\in [t_0, \Lambda t_0]$, where $M$ and $\nu$ are the
constants in \eqref{LTsectorial}, $\veps$ is the precision in the
evaluations of the Laplace transform $F$ and the elementary
operations in~\eqref{numinvLT},
$$
\Pi=\frac 1{\pi}
\sqrt{\frac{1+\sin(\alpha+d)}{(1-\sin(\alpha+d))^{2\nu-1}}}
$$
and
$$
Q = \max\{ 4L(\lambda t_0 \sin(\alpha-d)), \tau+L(\lambda t_0 \sin
\alpha) \},
$$
with $L(x)=1-\ln(1-e^{-x})$.
\end{theorem}

In case we have some reliable information about the errors in the
computation of the matrices $(zI-hA)^{-1}$, an $\veps$-depending
selection of $\tau$ and $\lambda$ improves the error bound
\eqref{cota} to $O(e^{-cK})$, i.e., a pure exponentially decaying
bound in $K$. In this situation, given $\veps$, $K$, and $\alpha, d$
fulfilling \eqref{angulos}, the parameters $\tau$ and $\mu$ are
given by
\begin{equation}\label{para_eps}
\tau = \frac {a(\theta^*)}{K}, \qquad \mu = \frac{2\pi d K
(1-\theta^*)}{\Lambda t_0 a(\theta^*)},
\end{equation}
where, for $\theta \in (0,1)$, $a(\theta)$ is the mapping
\begin{equation}\label{atheta}
a(\theta) = \acosh\left(\frac{\Lambda}{(1-\theta)\sin
\alpha}\right),
\end{equation}
and
\begin{equation}\label{theta_opt}
\theta^* = \min_{\theta \in (0,1)} \left( \veps \, e^{2\pi d K
(1-\theta)/a(\theta)} + e^{-2\pi d K \theta/a(\theta)} \right).
\end{equation}
Given $K \ge 1$, the expression to be minimize in \eqref{theta_opt}
represents the main part of the error bound obtained in
\cite{LoPSch} for any fixed $\theta \in (0,1)$. By choosing for
every $K$ the optimal value $\theta^*$ in \eqref{theta_opt}, an
error bound like $O(e^{-cK})$ is attained. We notice that the error
bound stated in Theorem~\ref{thm:err} is obtained for
$\theta^*=1-1/K$ in \eqref{para_eps}, instead of $\theta^*$ in
\eqref{theta_opt}. We refer to \cite{LoPSch} for details. There, the
case $\nu < 1$ in \eqref{LTsectorial} is also studied.


\section{Evaluation of the vector-valued mappings}\label{sec_eval}
In this section we apply some Laplace transformation formulas to
obtain a suitable representation of the operators $\phi_j(k,hA)$,
$\vfi_j(hA)$, and $\vfi_j(c_lhA)$ required in \eqref{multistep} and
\eqref{rk}.

Let us denote the Laplace transform of a mapping $f(\sigma)$ by
$F(z)=\cL[f](z)$, and the inverse Laplace transform by $f(\sigma) =
\cLinv[F](\sigma)$.

\subsection{Evaluation of the mappings required by the multistep
methods}\label{susec_mseval}

For $\phi_j$ in \eqref{phis_multistep} with $1 \le j\le k$, it holds
$$
\phi_j(k,\lambda) = \int_0^k e^{(k-\sigma)\lambda}
\left(\begin{array}{c} \sigma\\j \end{array} \right)\,d\sigma =
\cLinv[ \cL[f_0(\cdot,\lambda)] \times \cL [f_j] ] (k),
$$
where, for $\sigma > 0$,
\begin{equation}\label{f1fj}
f_0(\sigma,\lambda) = e^{\sigma \lambda} \qquad \mbox{and}\qquad
f_j(\sigma) = \left(\begin{array}{c} \sigma \\j \end{array}\right).
\end{equation}
For every $j \ge 1$ and $z\in \bC$, we define
\begin{equation}\label{LTPhi_j}
\Phi_j(z,\lambda) = \cL[f_0(\cdot,\lambda)](z) \times \cL [f_j](z) =
\frac{1}{z-\lambda}\times \cL [f_j](z).
\end{equation}
Then, for every $\lambda \in \bC$ and $j\ge 1$,
\begin{equation}\label{phisLinv}
\phi_j(k,\lambda) = \cLinv [\Phi_j(\cdot,\lambda)](k).
\end{equation}
For $j=0$
$$
\phi_0(k,\lambda) = e^{k\lambda} = \cLinv \left(\frac {1}{\cdot
-\lambda}\right)(k),
$$
and thus we define
\begin{equation}\label{Phicero}
\Phi_0(z,\lambda) =\frac {1}{z-\lambda}.
\end{equation}
For $\lambda$ scalar, the mappings $\Phi_j(z,\lambda)$, with $1\le j
\le 4$ are given by
\begin{equation}\label{scalarformulae}
\begin{array}{ll}
\Phi_1(z,\lambda) = \displaystyle\frac{1}{z(z-\lambda)}, \quad  &
\Phi_2(z,\lambda) = \displaystyle \frac{1}{z^2(z-\lambda)}, \\[10pt]
\Phi_3(z,\lambda) = \displaystyle \frac{2-z}{2z^3(z-\lambda)}, \quad
& \Phi_4(z,\lambda) = \displaystyle
\frac{3-3z+z^2}{3z^4(z-\lambda)}.
\end{array}
\end{equation}

In order to evaluate $\phi_j(k,hA)$, $0\le j\le 4$, we propose to
use the formulas in \eqref{scalarformulae} with $hA$ instead of
$\lambda$ and perform the inversion of the Laplace transform to
approximate the original mappings at $\sigma=k$. In this way, the
Laplace transforms we need to invert are:
\begin{equation}\label{operator_multistep}
\begin{array}{lcl}
\Phi_0(z,hA) &=& (zI-hA)^{-1}, \\[1em]
\Phi_1(z,hA) &=& \displaystyle \frac 1z (zI-hA)^{-1}, \\[1em]
\Phi_2(z,hA) &=& \displaystyle \frac 1{z^2}(zI-hA)^{-1},\\[1em]
\Phi_3(z,hA) &=& \displaystyle \frac{2-z}{2z^3}(zI-hA)^{-1},\\[1em]
\Phi_4(z,hA) &=& \displaystyle \frac{3-3z+z^2}{3z^4}(zI-hA)^{-1}.
\end{array}
\end{equation}
Although the formulas in \eqref{operator_multistep} are derived just
formally, we notice that they can be justified by combining the
Cauchy integral formula with the inversion formula for the Laplace
transform. More precisely, for suitable contours $\Gamma_1$ and
$\Gamma_2$ in the complex plane, both laying in the resolvent set of
$A$, it holds
\begin{equation}\label{proofformula}
\begin{array}{lcl}
\phi_j(k,hA) &=& \displaystyle \frac{1}{2\pi i} \int_{\Gamma_1} \phi_j(k,\lambda) (\lambda I-hA)^{-1}\, d\lambda \\[10pt]
&=&  \displaystyle  \frac{1}{2\pi i} \int_{\Gamma_1} \bigg( \frac{1}{2\pi i} \int_{\Gamma_2} e^{k \xi} \Phi_j(\xi,\lambda) \, d\xi \bigg) (\lambda I-hA)^{-1} \,d\lambda \\[10pt]
&=& \displaystyle  \frac{1}{2\pi i} \int_{\Gamma_2} e^{k\xi} \bigg(\frac{1}{2\pi i} \int_{\Gamma_1} \Phi_j(\xi,\lambda) (\lambda I-hA)^{-1}\, d\lambda \bigg) \, d\xi\\[10pt]
&=& \displaystyle  \frac{1}{2\pi i} \int_{\Gamma_2} e^{k\xi}
\Phi_j(\xi,hA) \, d\xi = \cLinv [\Phi_j(\cdot,hA)](k).
\end{array}
\end{equation}

Due to \eqref{sectorial}, all the Laplace transforms $\Phi_j(k,hA)$
in \eqref{operator_multistep} are sectorial, since they satisfy
\eqref{LTsectorial} with $\gamma^* = \max\{0,\gamma\}$ and
$\delta^*= \delta$, for $\gamma$ and $\delta$ in \eqref{sectorial}.
We notice that, for all $j$, the resulting bounds in
\eqref{LTsectorial} are independent of $h$.

Thus, we can compute the operators $\phi_j(k,hA)$, $0\le j\le k$, by
using the method described in Section~\ref{sec_LT} to compute the
inverse Laplace transforms of the mappings $\Phi_j(z,hA)$ in
\eqref{LTPhi_j}. We notice that the inverse Laplace transforms need
to be approximated only at the fix value $\sigma=k$, which is
specially favorable for the application of the inversion method (see
the bound in Theorem~\ref{thm:err}). Then, we set $\Lambda = 1$,
$t_0=k$, and select the parameters $\tau$ and $\mu$ following
Theorem~\ref{thm:err}. The selection of $\alpha$ and $d$ is more
heuristic and a good choice is $\alpha \approx \frac 12(\frac
{\pi}2-\delta)$ and $d$ slighly smaller than $\alpha$. For example,
if $\delta = 0$ in \eqref{LTsectorial}, good values are around
$\alpha=0.7$ and $d=0.6$. Next, we compute the quadrature weights
$w_{\ell}$ and nodes $z_{\ell}$ in \eqref{pesos_nodos} and
approximate the operators in \eqref{multistep} by
\begin{equation}\label{aproxphis}
\phi_j(k,hA) \approx \sum_{\ell=-K}^{K} w_{\ell} e^{k z_{\ell}}
\Phi_j(z_{\ell},hA).
\end{equation}
The sum in \eqref{aproxphis} can be halved in case of symmetry like
in \eqref{numinvLT_real}.

As we already mentioned in the Introduction, the computation of all
the required operators $\phi_j(k,hA)$, $0\le j\le k$, can be carried
out before the time stepping of the exponential method begins. Thus,
if we use the method of lines and apply the exponential method to
some spatial discretization of \eqref{ivp}, only the matrix-vector
products in \eqref{multistep} need to be computed at every time
step.

\subsection{Evaluation of the mappings required by the Runge--Kutta
methods}\label{subsec_rkeval} For $\vfi_j$ in \eqref{phis_rk}, $j\ge
1$ and $t>0$, we have
$$
\vfi_j(\lambda) = \int_0^1 e^{(1-\sigma)\lambda}
\frac{\sigma^{j-1}}{(j-1)!}\,d\sigma = \cLinv[
\cL[g_0(\cdot,\lambda)] \times \cL [g_j] ] (1),
$$
where, for $\sigma > 0$ and $\lambda \in \bC$,
\begin{equation}\label{g1gj}
g_0(\sigma,\lambda) = e^{\sigma \lambda} \qquad \mbox{and}\qquad
g_j(\sigma) = \frac{\sigma^{j-1}}{(j-1)!}.
\end{equation}
For every $j \ge 1$ and $z\in \bC$, we define
\begin{equation}\label{LTPsi_j}
\Psi_j(z,\lambda) = \cL[g_0(\cdot,\lambda)](z) \times \cL [g_j](z) =
\frac{1}{z^j(z-\lambda)}
\end{equation}
and $\Psi_0(z,\lambda)=(z-\lambda)^{-1}$. Then, for every $\lambda
\in \bC$ and $j\ge 0$,
\begin{equation}\label{vfisLinv}
\vfi_j(\lambda) = \cLinv [\Psi_j(\cdot,\lambda)](1).
\end{equation}

The same argument as in \eqref{proofformula} justifies the
computation of the operators $\vfi_j(hA)$ and $\vfi_j(c_lhA)$, $j\ge
0$, $2\le l\le s$, by performing the inversion of the Laplace
transforms
\begin{equation}\label{operator_rk}
\Psi_j(z,\beta hA) = \displaystyle \frac 1{z^j} (zI- \beta hA)^{-1},
\qquad j \ge 0,\quad \beta =1,\ c_l,
\end{equation}
to approximate the original mappings at $\sigma=1$. If
\eqref{sectorial} holds, the Laplace transforms in
\eqref{operator_rk} are also the sectorial in the sense of
\eqref{LTsectorial} and we can use the inversion method of
\cite{LoPSch}.

As in the case of the methods in \eqref{multistep}, the computation
of all the required operators $\vfi_j(hA)$ and $\vfi_j(c_lhA)$,
$j\ge 0$, can be carried out before the time stepping of the
exponential method begins.

\begin{nota}
In general, we can always evaluate a mapping $\phi(hA)$ of the form
of \eqref{Phigeneral} by using the numerical inversion of the
Laplace transform, just by noticing that $\phi(hA)$ is the inverse
Laplace transform at $\sigma = n$ of a mapping $\Phi(z,hA)$ like in
\eqref{LTPhi_general},
$$
\Phi(z,hA)=R(z)(zI-hA)^{-1},
$$
with $R(z) = \cL[p](z)$, a scalar rational function of $z$.
\end{nota}

The above Remark implies that our algorithm can be used to implement
other kinds of exponential methods, different than those in
\cite{CalPal06,HoOs05exp}, as long as they require the evaluation of
mappings of the form of \eqref{Phiexponencial} and
\eqref{Phigeneral}.

\section{Evaluation of the scalar mappings}\label{sec_evalscalar}
As we already mentioned in the Introduction, we can also apply the
inversion of the Laplace transform to evaluate with accuracy the
scalar mappings $\phi(\lambda)$ in \eqref{Phigeneral}. In this
section we consider with some detail the evaluation of the mappings
\begin{equation}\label{gj}
g_j(m,\lambda) = \int_0^{m} e^{(m-\sigma)\lambda}\sigma^{j-1}
\,d\sigma, \qquad
 j\ge 1, \quad  m\in \bN,
\end{equation}
by means of the quadrature formula \eqref{quadscalar}. The result
provided by \eqref{quadscalar} does not depend on the size of
$\lambda$, but the formula is only useful in principle for values of
$\lambda$ inside a sector of the form $|\arg (\lambda-\gamma)| \ge
\pi-\delta$. However, using that
\begin{equation}\label{g1_minus}
e^{m\lambda}g_1(m,-\lambda) = g_1(m,\lambda), \qquad \lambda \in
\bC,
\end{equation}
and
\begin{equation}\label{recursion_gj}
g_{j+1}(m,\lambda) = \frac{jg_{j}(m,\lambda)-m^{j}}{\lambda}, \qquad
j\ge 1,\quad m\in \bN,
\end{equation}
it is easy to see by induction that, for $m\in \bN$ and $\lambda\in
\bC$,
\begin{equation}\label{gj_minus}
e^{m\lambda}g_{j}(m,-\lambda)=\sum_{\ell=1}^{j}
\left(\begin{array}{c} j-1
\\ \ell-1
\end{array} \right) (-1)^{\ell-1} m^{j-\ell}g_{\ell}(m,\lambda),\quad j\ge 1.
\end{equation}
Thus, we can compute
\begin{equation}\label{gj_minus_LT}
g_j(m,-\lambda)= \displaystyle e^{-m\lambda} \cL^{-1} \left[
G^*_{j}(\cdot,\lambda)\right](m),
\end{equation}
with
\begin{equation}\label{LTgj_minus}
G^*_j(z,\lambda) = \frac{1}{z^{j}(z-\lambda)}\sum_{\ell=1}^{j}
\left(\begin{array}{c} j-1
\\ \ell-1
\end{array} \right)\ell! (-1)^{\ell-1} (m z)^{j-\ell}, \quad j\ge 1.
\end{equation}
which provides a stable formula to approximate $g_j(m,-\lambda)$ for
$\lambda$ inside a proper sector $|\arg(\lambda-\gamma)|>\pi-\delta$
and moderate size.

In Table~\ref{tab1} we show the error obtained in the evaluation of
$\vfi(\lambda)=g_1(1,\lambda)$ in \eqref{vfi_1} for different values
of $\lambda$ in the interval $[-1,1]$. For $\lambda < 0$, we applied
the inversion formula \eqref{numinvLT} with $t=1$ and
\begin{equation}\label{LTg1}
F(z)=G_1(z,\lambda) = \frac1{z(z-\lambda)},
\end{equation}
which, for these values of $\lambda$, fulfils \eqref{LTsectorial}
with $\delta = 0$, $\gamma=0$, and $\nu=2$. We assumed that the
evaluations of $G_1$ can be carried out in MATLAB up to machine
accuracy and thus we set $\veps = 2.2204 \times 10^{-16}$. Then, we
computed the quadrature weights and nodes in \eqref{pesos_nodos}
following \eqref{para_eps}--\eqref{theta_opt} with $\Lambda=1$.
Setting $\alpha = 0.7$ and $d=0.6$, we obtained $\theta = 0.693$,
for $K=15$, and $\theta=0.793$, for $K=25$. In Table~\ref{tab1} we
can see that $K=25$ is enough to attain almost the machine accuracy
of MATLAB in the evaluations of $\vfi(\lambda)$. For positive values
of $\lambda$, we used \eqref{g1_minus} with $m=1$.

\begin{table}[t]
\caption{Computation of $\vfi(\lambda)$ in \eqref{vfi_1} for
$\lambda\in[-1,1]$ by using formulas \eqref{numinvLT} and
\eqref{g1_minus}. We show the absolute error obtained in MATLAB with
$K=15$ and $K=25$.}\label{tab1}
\begin{center}
\begin{tabular}{|c|c|c|c|c|c|}
\hline $\lambda < 0$ & $K=15$ & $K=25$ & $-\lambda$ & $K=15$ & $K=25$ \\
\hline
-1 & 1.5050e-12 & 1.3323e-15 & 1 & 1.5050e-12 & 3.3307e-15 \\
 -1e-1 & 1.5227e-12 & 3.2196e-15 & 1e-1 & 1.5227e-12 & 3.5527e-15 \\
 -1e-2 & 1.4243e-12 & 4.4409e-15 & 1e-2 & 1.4243e-12 & 4.6629e-15 \\
 -1e-3 & 1.3750e-12 &  1.3323e-15 & 1e-3 & 1.3750e-12 & 1.3323e-15 \\
 -1e-4 & 1.3738e-12 & 1.7764e-15 & 1e-4 & 1.3738e-12 & 1.7764e-15 \\
 -1e-5 & 1.3747e-12 & 3.6637e-15 & 1e-5 & 1.3747e-12 & 3.7748e-15 \\
 -1e-6 & 1.3748e-12 & 3.6637e-15 & 1e-6 & 1.3748e-12 & 3.7748e-15 \\
 -1e-7 & 1.3695e-12 & 1.9984e-15 & 1e-7 & 1.3695e-12 & 1.9984e-15 \\
 -1e-8 & 1.3717e-12 & 1.1102e-16 & 1e-8 & 1.3717e-12 & 2.2204e-16 \\
 -1e-9 & 1.3715e-12 & 1.1102e-16 & 1e-9 & 1.3715e-12 & 0\\
 -1e-10 & 1.3711e-12 &       0 & 1e-10 & 1.3711e-12 & 0 \\
 -1e-11 & 1.3711e-12 &     0 & 1e-11 & 1.3711e-12 & 0 \\
 -1e-12 & 1.3715e-12 & 1.1102e-16 & 1e-12 & 1.3715e-12 & 0 \\
 -1e-13 & 1.3712e-12 &        0 & 1e-13 & 1.3712e-12 & 2.2204e-16 \\
\hline
\end{tabular}
\end{center}
\end{table}

\section{Numerical illustrations}\label{sec_experiments}
In this section we test our algorithm by considering the same
examples as in \cite{CalPal06} and \cite{HoOs05exp}.

\subsection{Example for the multistep exponential methods}
Our first example is the problem considered in \cite{CalPal06}
\begin{equation}\label{ex1_cp}
u_t(x,t) = u_{xx}(x,t) + \left(\int_0^1 u(s,t)\,ds \right) u_x(x,t)
+ g(x,t),
\end{equation}
for $x\in [0,1]$ and $t\in[0,1]$, subject to homogeneous Dirichlet
boundary conditions and with $g(x,t)$ such that the exact solution
to \eqref{ex1_cp} is $u(x,t)=x(1-x)e^t$

The spatial discretization of \eqref{ex1_cp} is carried out by using
standard finite differences with $J = 512$ spatial nodes, centered
for the approximation of $u_x$. The nonlocal term is approximated by
means of the composite Simpson's formula.

To integrate in time the semidiscrete problem we use
\eqref{multistep} with $k=1, 2, 3$ and $4$, so that $A$ is the
$(J-1)\times(J-1)$ matrix
$$
A = J^2 {\rm tridiag\,}([1,-2,1]).
$$
We approximate the matrices $\phi_j(k,hA)$, $0\le j \le k$, required
in \eqref{multistep} by applying the quadrature rule
\eqref{aproxphis}. To avoid an extra source of error, the initial
values $u_1,\dots,u_{k-1}$ are computed from the exact solution. In
a less academic example, these values can be computed by means of a
one-step method of sufficiently high order or by performing the fix
point iteration proposed in \cite{CalPal06}.

\begin{figure}
\centering
\includegraphics[width=0.49\textwidth]{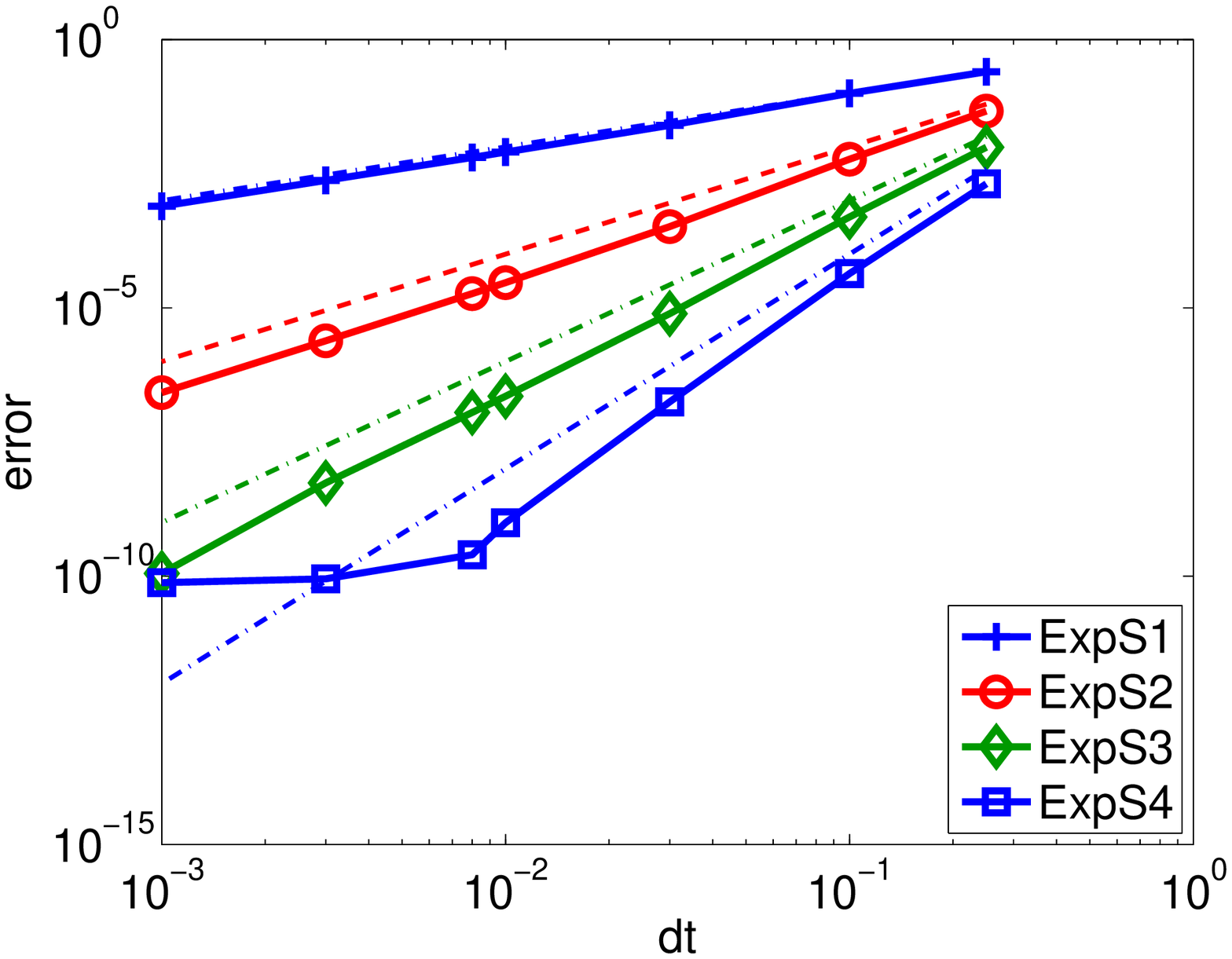}
\includegraphics[width=0.49\textwidth]{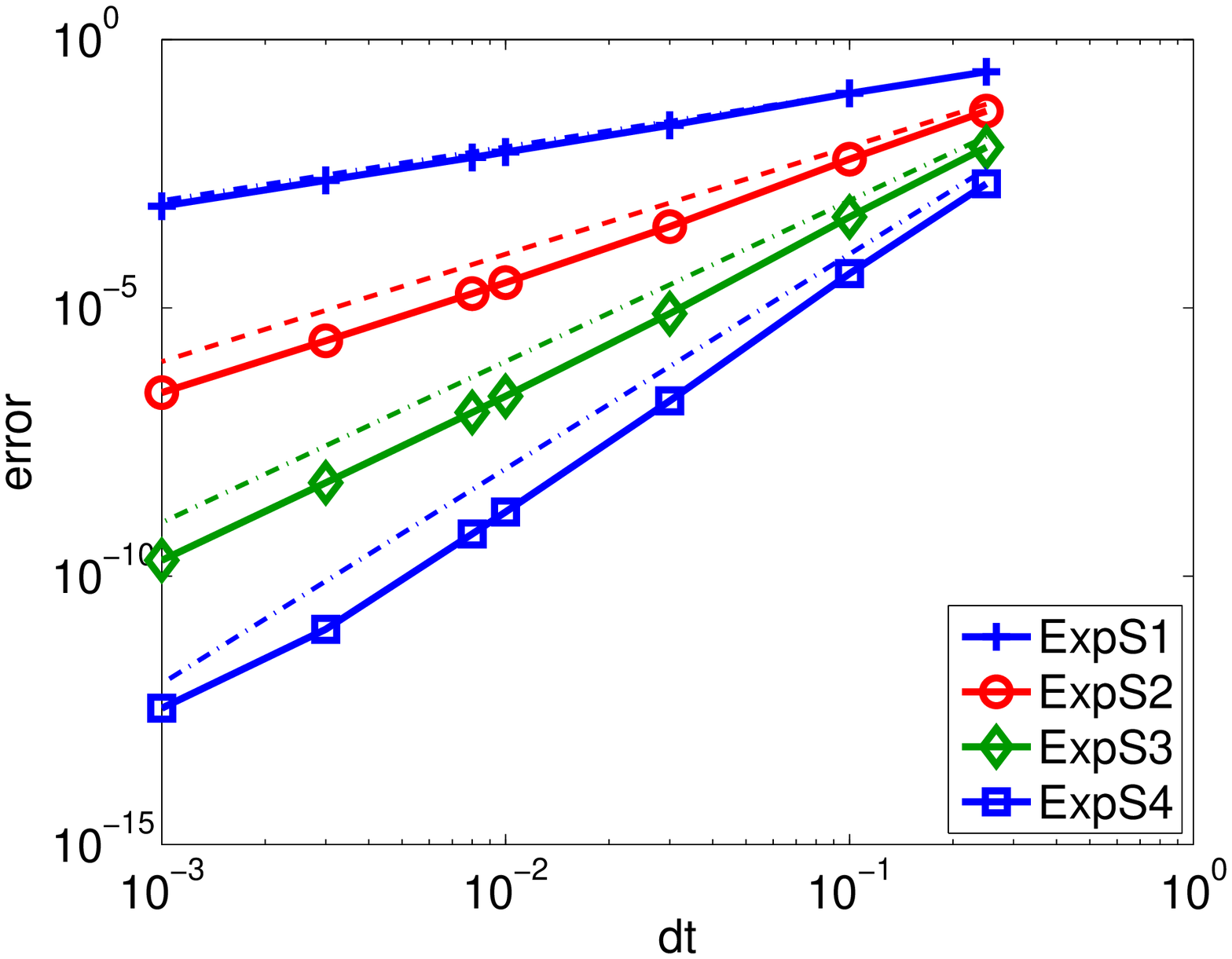}
\caption{Error of exponential multistep methods \eqref{multistep}
applied to \eqref{ex1_cp}, for $k=1, 2, 3$, and $4$. Left: With
$K=25$ quadrature nodes on the hyperbolas, Right: With $K=35$.}
\label{fig:ex1_cp}
\end{figure}

In Figure~\ref{fig:ex1_cp} we show the error versus the stepsize at
$t=1$, measured in a discrete version of the norm $\|\cdot\|_{1/2}$,
for $K=25$ and $K=35$ in \eqref{aproxphis}. We see that for $K=35$
the full precision is achieved for all the methods implemented;
cf.~\cite[Section 6]{CalPal06}. In Figure~\ref{fig:ex1_cp} we also
show lines of slope 1, 2, 3 and 4, to visualize the order of
convergence. In fact, we can observe that the order of convergence
for this example is slightly higher than the one predicted in
\cite{CalPal06}, approximately 2.15, 3.15 and 4.15 for $k=2, 3$ and
$4$, respectively, instead of sharp order $k$. A further study of
this behavior is beyond the scope of this paper.

\subsection{Examples for the exponential Runge--Kutta methods}
For the following two examples we consider the problems and some of
the integrators proposed in \cite{HoOs05exp}. Following the notation
in \cite{HoOs05exp}, in the Butcher tableaus below we use the
abbreviations
\begin{equation}\label{abbrev}
\vfi_i = \vfi_i(hA), \quad \mbox{and }\quad  \vfi_{i,j} = \vfi_{i,j}(hA) = \vfi_i(c_jhA), \qquad 2\le j \le s.
\end{equation}

Thus, our second example is
\begin{equation}\label{ex1_ho}
u_t(x,t) = u_{xx}(x,t) + \frac 1{1 + u(x,t)^2} + g(x,t),
\end{equation}
for $x\in[0,1]$ and $t\in [0,1]$, subject to homogeneous Dirichlet
boundary conditions and with $g(x,t)$ such that the exact solution
to \eqref{ex1_ho} is again $u(x,t)=x(1-x)e^t$.

We discretize this problem in space by standard finite differences
with $J=200$ grid points. For the time integration of the
semidiscrete problem, we implemented \eqref{rk} with $s=1$, the
second-order method
\begin{equation}\label{rk2}
\begin{array}{r|cc}
0 &  &  \\
\frac 12 & \frac 12 \vfi_{1,2} & \\[.3em] \hline
 & 0 & \vfi_1
\end{array}
\end{equation}
the third-order method
\begin{equation}\label{rk3}
\begin{array}{r|ccc}
0 &  &  &\\[.3em]
\frac 13 & \frac 13 \vfi_{1,2} & & \\[.3em]
\frac 23 & \frac 23 \vfi_{1,3} - \frac 43 \vfi_{2,3} & \frac 43\vfi_{2,3} \\[.3em] \hline
& & & \\[-1em]
 & \vfi_1 - \frac 32\vfi_2 & 0 & \frac 32 \vfi_2
\end{array}
\end{equation}
and the fourth-order one
\begin{equation}\label{rk4}
\begin{array}{r|ccccc}
0 &  &  & & &\\[.3em]
\frac 12 & \frac 12 \vfi_{1,2} & & & &\\[.3em]
\frac 12 & \frac 12 \vfi_{1,3} - \vfi_{2,3} & \vfi_{2,3} & & &\\[.3em]
1 & \vfi_{1,4} - 2 \vfi_{2,4} & \vfi_{2,4} & \vfi_{2,4} & &\\[.3em]
\frac 12 & \frac 12 \vfi_{1,5}-2a_{5,2}-a_{5,4} & a_{5,2} & a_{5,2} & a_{5,4} &
 \\[.3em] \hline
& & & \\[-1em]
& \vfi_1-3\vfi_2+4\vfi_3 & 0 & 0 & -\vfi_2 + 4\vfi_3 & 4\vfi_2-8\vfi_3
\end{array}
\end{equation}
with
$$
a_{5,2} = \frac 12 \vfi_{2,5} - \vfi_{3,4} + \frac 14 \vfi_{2,4} - \frac 12 \vfi_{3,5}
$$
and
$$
a_{5,4} = \frac 14 \vfi_{2,5}-a_{5,2}.
$$

For the implementation of \eqref{rk2} we need to invert four
different Laplace transforms of the form of \eqref{operator_rk}, to
approximate $\vfi_0 \left(\frac h2 A \right), \vfi_0(hA)$,
$\vfi_1\left(\frac h2 A \right)$, and $\vfi_1(hA)$. The
implementation of both \eqref{rk3} and \eqref{rk4} requires the
inversion of eight Laplace transforms.

\begin{figure}
  \centering
    \includegraphics[width=0.49\textwidth]{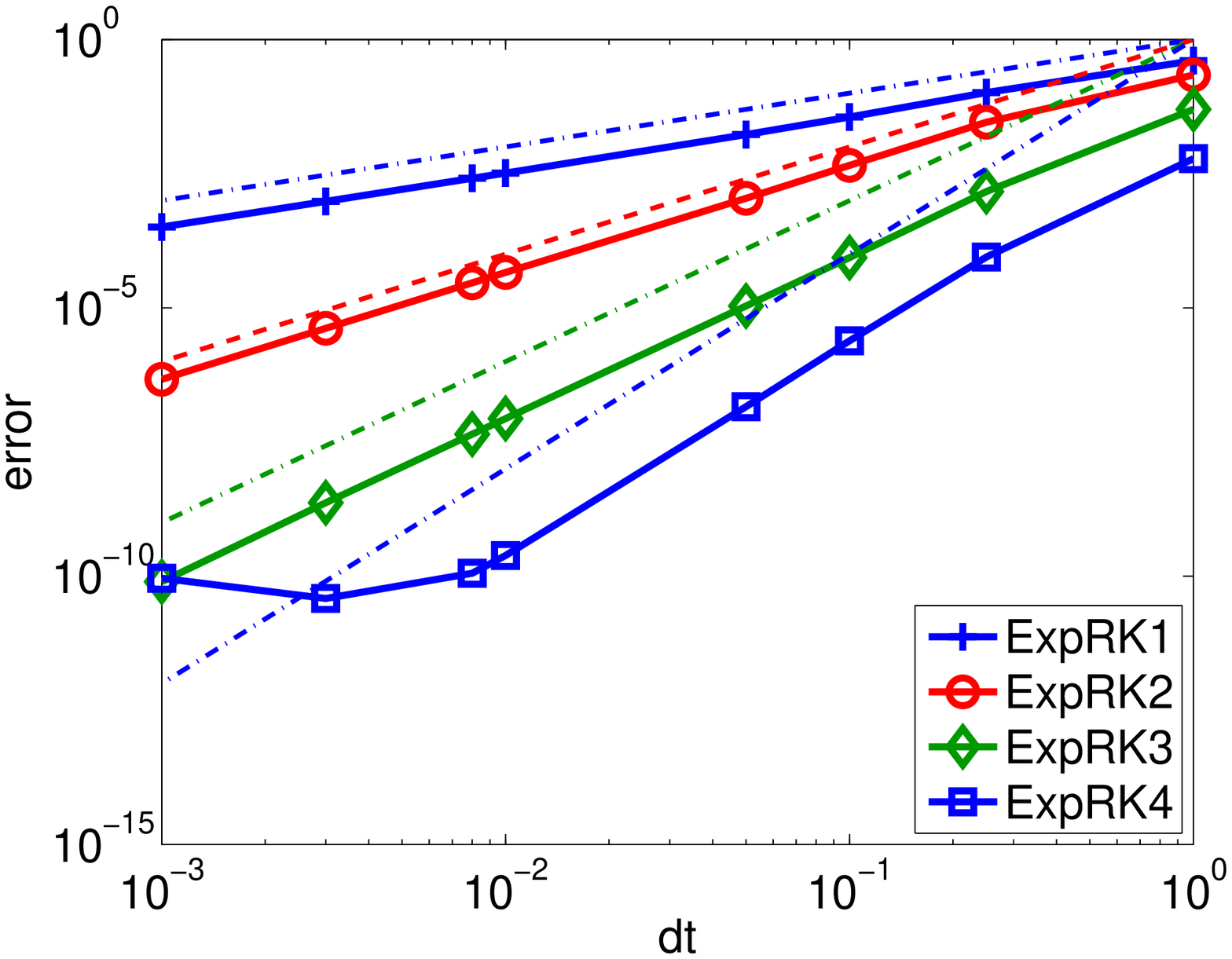}
    \includegraphics[width=0.49\textwidth]{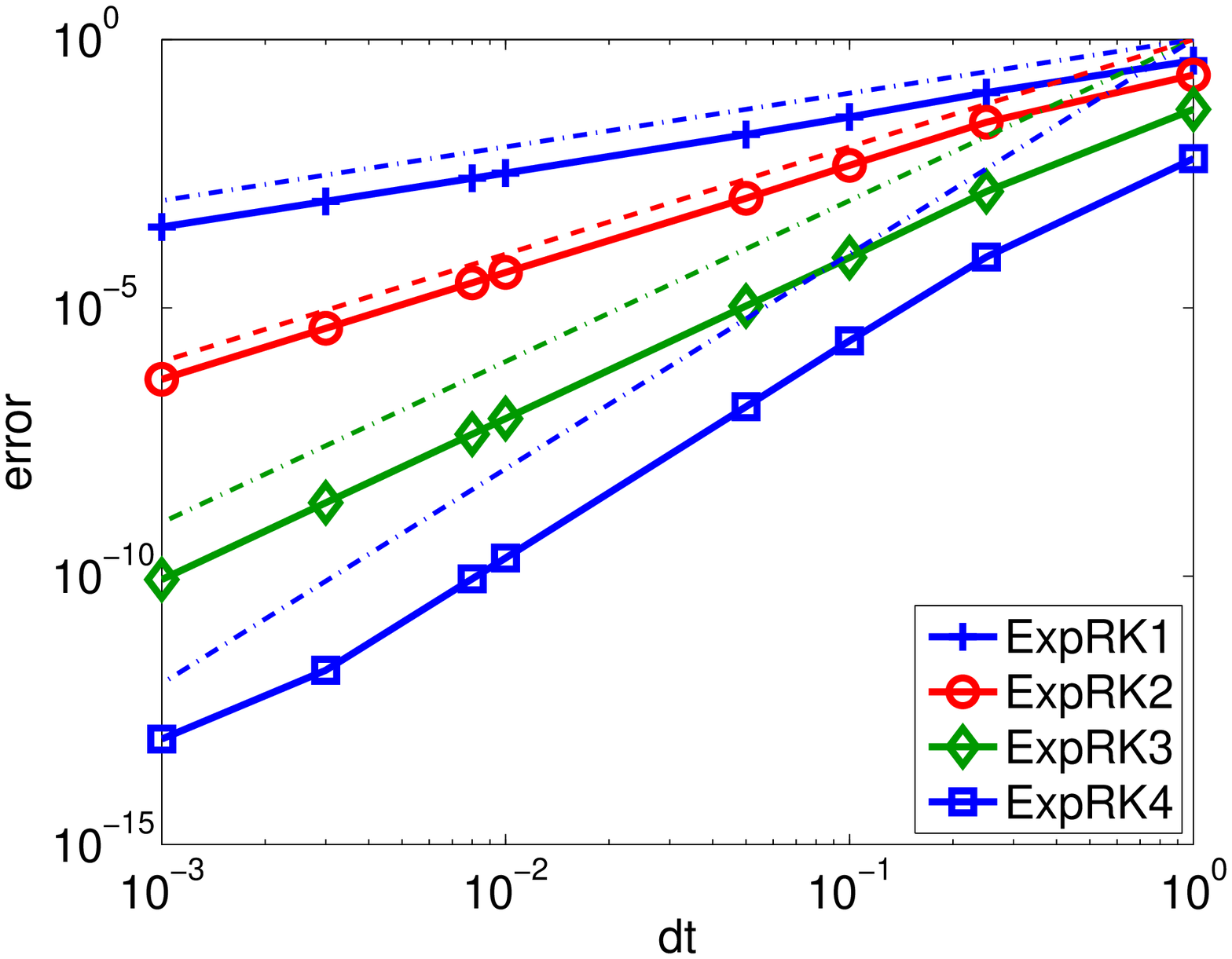}
  \caption{Error of Runge--Kutta methods \eqref{multistep} with $s=1$, \eqref{rk2}, \eqref{rk3}, and \eqref{rk4},
   applied to \eqref{ex1_ho}. Left: With $K=25$ quadrature nodes in \eqref{aproxphis}, Right: With $K=35$.}
  \label{fig:ex1_ho}
\end{figure}

In Figure~\ref{fig:ex1_ho} we show the error at $t=1$ versus the
stepsize, measured in the maximum norm. The expected order of
convergence for this example is $k$ for the $k$-order method. In
order to check our algorithm, we added lines with the corresponding
slopes in Figure~\ref{fig:ex1_ho}. We can see that also for this
kind of methods we attain full precision for $K=35$ in
\eqref{aproxphis}.

Finally we consider the second example in \cite{HoOs05exp}
\begin{equation}\label{ex2_ho}
u_t(x,t) = u_{xx}(x,t) + \int_0^1 u(s,t)\,ds + g(x,t),
\end{equation}
for $x\in[0,1]$ and $t\in [0,1]$, subject also to homogeneous
Dirichlet boundary conditions and with $g(x,t)$ such that the exact
solution to \eqref{ex2_ho} is $u(x,t)=x(1-x)e^t$.

\begin{figure}
  \centering
    \includegraphics[width=0.49\textwidth]{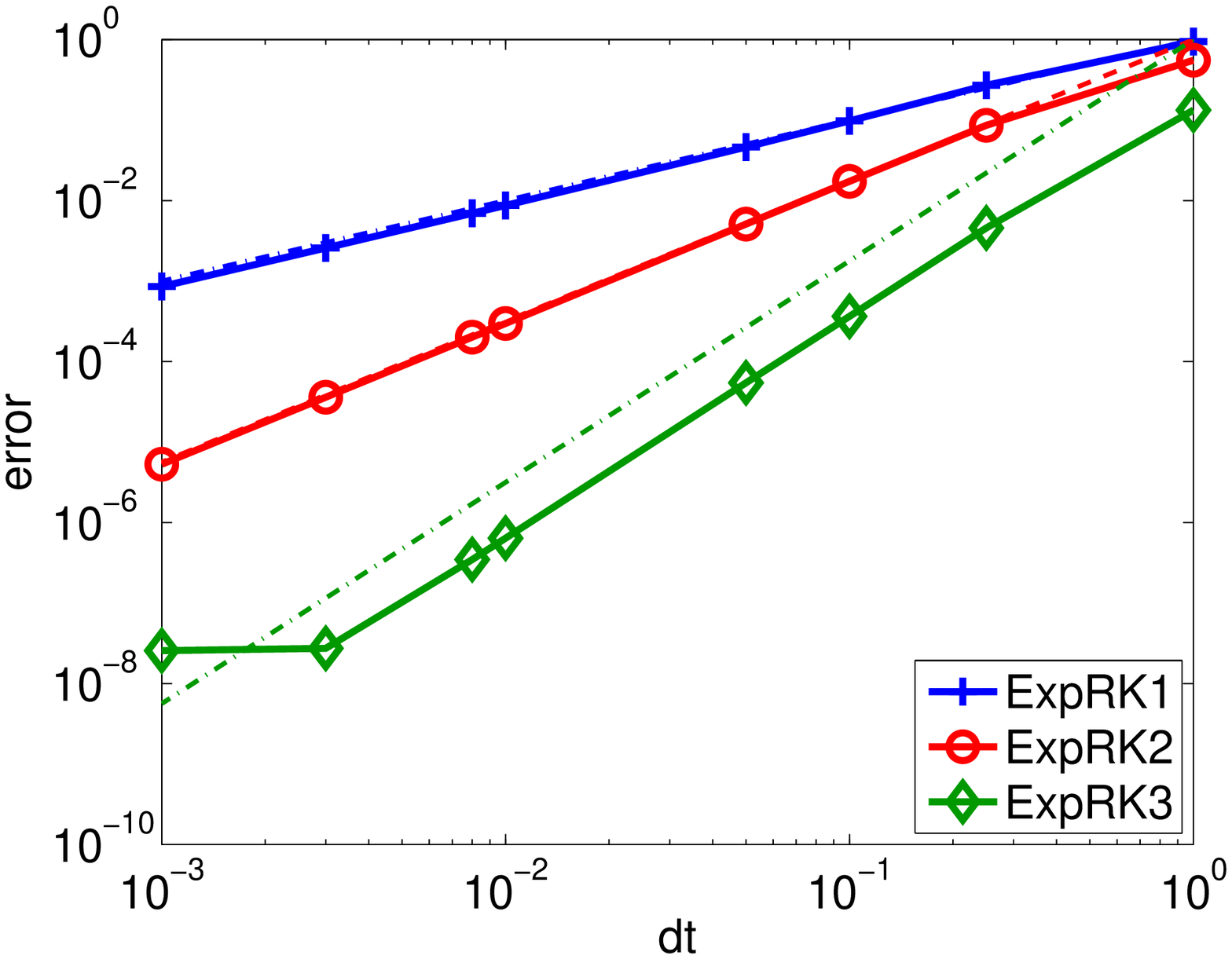}
    \includegraphics[width=0.49\textwidth]{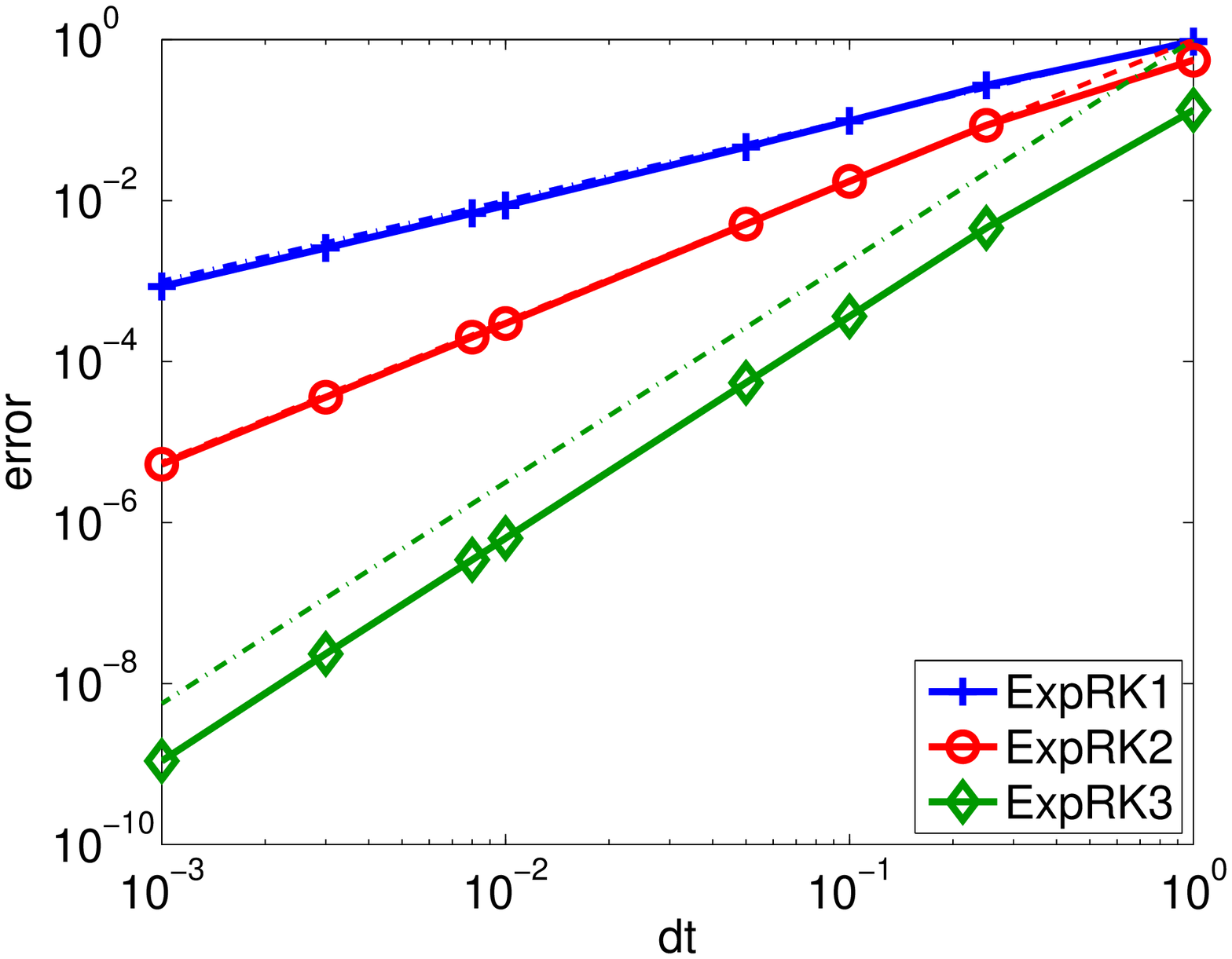}
  \caption{Error of Runge--Kutta methods \eqref{rk} with $s=1$, \eqref{rk2}, and \eqref{rk3},
   applied to \eqref{ex2_ho}. Left: With $K=20$ quadrature nodes in \eqref{aproxphis}, Right: With $K=25$.}
  \label{fig:ex2_ho}
\end{figure}

We discretize this problem in space as in the previous example
\eqref{ex1_ho}, and use the composite Simpson's rule for the
approximation of the nonlocal term. For the time integration, we use
\eqref{rk} with $s=1$, \eqref{rk2}, and \eqref{rk3}. In
Figure~\ref{fig:ex2_ho} we show the error at $t=1$, measured in a
discrete version of the norm $\|\cdot\|_{1/2}$;~cf.~\cite[Section
6]{HoOs05exp}. The expected order of convergence for this example is
1 for the first order method, $1.75$ for \eqref{rk2}, and 2.75 for
\eqref{rk3}. In order to test our algorithm for the exponential
methods, we added lines with the corresponding slopes in
Figure~\ref{fig:ex2_ho}. We refer to \cite{HoOs05exp} for a detailed
explanation of the order reduction phenomenon in this example. The
less stringent accuracy requirements in this case allow us to
achieve full precision with only $K=25$ quadrature nodes in
\eqref{aproxphis}.

\section{Conclusions}
In this paper we derived a way to approximate the exponential-like
operators required for the implementation of different kinds of
exponential me\-thods present in the recent literature for the time
integration of semilinear problems. The approach is based on the
numerical inversion of the Laplace transform and its application is
restricted to parabolic problems. When applicable, the proposed
algorithm is shown to be very efficient, both at the scalar level
and the operator level, and it is quite simple to implement. Apart
from the error bound stated in Theorem~\ref{thm:err}, which is a
partial result of those proved in \cite{LoPSch}, we tested the
algorithm with several academic examples from \cite{CalPal06} and
\cite{HoOs05exp} and with the evaluation of the prototypical mapping
$\vfi$ in \eqref{vfi_1}.

The comparison of our algorithm with other techniques present in the
literature (see the Introduction) is beyond the scope of the present
paper. Some testing for comparison of CPU times and storage
requirements with problems of different complexity could be the
subject of future work.

\section{Acknowledgements} The author is grateful to Enrique Zuazua
for proposing a good problem and for useful discussions during the
preparation of the manuscript. She also thanks C\'esar Palencia for
suggesting the use of exponential methods and Christian Lubich for
his valuable comments on a preliminary version of the paper.

\bibliographystyle{amsplain}

\end{document}